\documentclass[12pt]{amsproc}%

\usepackage{amssymb}
\usepackage{amsmath}
\usepackage{amsthm}
\usepackage{amscd}
\usepackage[margin=1in]{geometry}%
\usepackage{hyperref}
\usepackage{appendix}

\theoremstyle{plain}
\newtheorem{theorem}{Theorem}[section]
\newtheorem{proposition}{Proposition}[section]
\newtheorem{corollary}{Corollary}[section]
\newtheorem{lemma}{Lemma}[section]

\theoremstyle{remark}

\newtheorem{remark}{Remark}[section]
\newtheorem{examples}{Examples}[section]

\DeclareMathOperator{\clos}{clos}
\DeclareMathOperator{\lin}{span}
\DeclareMathOperator{\diverg}{div}

\DeclareMathOperator{\curl}{curl}

\begin{document}

\title{Vector analysis on fractals and applications}
\author[Hinz]{Michael Hinz$^1$}
\address{Mathematisches Institut, Friedrich-Schiller-Universit\"at Jena, Ernst-Abbe-Platz 2, 07737, Germany and Department of Mathematics,
University of Connecticut,
Storrs, CT 06269-3009 USA}
\email{Michael.Hinz.1@uni-jena.de and Michael.Hinz@uconn.edu}
\thanks{$^1$Research supported in part by NSF grant DMS-0505622 and by the Alexander von Humboldt Foundation Feodor (Lynen Research Fellowship Program)}
\author[Teplyaev]{Alexander Teplyaev$^2$}
\address{Department of Mathematics, University of Connecticut, Storrs, CT 06269-3009 USA}
\email{Alexander.Teplyaev@uconn.edu}
\thanks{$^2$Research supported in part by NSF grant DMS-0505622}
\date{\today}

\begin{abstract}
The paper surveys some recent results concerning vector analysis on fractals. We start with a local regular Dirichlet form and use the framework of $1$-forms and derivations introduced by Cipriani and Sauvageot to set up some elements of a related vector analysis in weak and non-local formulation.  This allows to study various scalar and vector valued linear and non-linear partial differential equations on fractals that had not been accessible before. Subsequently a stronger (localized, pointwise or fiberwise) version of this vector analysis can be developed, which is related to previous work of Kusuoka, Kigami, Eberle, Strichartz, Hino, Ionescu, Rogers, R\"ockner, and the authors.  
\tableofcontents
\end{abstract}
\maketitle

\section{Introduction}

In the present article we survey some recent results concerning vector analysis based on symmetric local regular Dirichlet forms on locally compact separable metric spaces, cf. \cite{BH91,FOT94,RocknerMa,RW}. 
The notions and results we discuss have been introduced in the papers \cite{HRT, HTa, HTb}. 
They are based on the approach to differential $1$-forms as proposed by Cipriani and Sauvageot in \cite{CS03} in much greater generality, and later investigated by several authors, \cite{CGIS11, CS09, H11, IRT}. 
The constructions are sufficiently robust to apply to symmetric diffusions on fractals such as 
p.c.f. self-similar sets \cite{Ki93, Ki01,Str06}, 
nested fractals \cite{Li}, 
finitely ramified fractals \cite{S97, T08}, 
generalized Sierpinski carpets \cite{BB99, KuZh92, BBKT}, 
spaces of Berlow-Evans-Laakso type \cite{S-POTA, ST-BE},  
and some random fractals \cite{Ham92, Ham97}. 
As they are based on Dirichlet forms they also apply to classical situations such as Euclidean spaces, domains with sufficiently regular boundary and smooth compact Riemannian manifolds. In these cases we recover well-known results. 

A general theme motivating our studies consists of the questions for which elements of differential geometry and vector analysis one can find analogs built solely upon the notion of energy and how these analogs can be used to formulate and study physical models on non-smooth spaces. The space $\mathcal H$ of $1$-forms as constructed in \cite{CS03, CS09} is  a Hilbert space. Therefore one can identify $1$-forms and vector fields, and furthermore introduce other notions of vector analysis, as recently done in \cite{HRT} (which generalizes earlier approaches to vector analysis on fractals, see \cite{Ki93,Ki08,Ku93,PT,St00,T08}). This is a part of a  comprehensive program to introduce and study  vector equations on general non-smooth spaces   which carry a diffusion process (or, equivalently, a local regular Dirichlet form). 

\def\arti{and references therein}

Much of the existing literature on analysis on fractals has been concerned with the primary problems of construction diffusions on fractals, 
\cite[\arti]{BB89, BBKT, BP88, G87, HMT, Ki93, KumSt05, Ku87, Ku89, KuZh92, Li, S97}, 
studying their heat kernel decay, 
\cite[\arti]{BB92, BK01, GT12, Ka12, Ka12b, Ki09, Kum93}, 
their potential theory, 
\cite[\arti]{BBK06, BGK11, B11, IRS, M91}, 
their spectral properties, 
\cite[\arti]{ADT1, ADT2, e1, e2, DGV, FSh92, Ham11, HSTZ, KiLa93, KL2, Ka10, La91, LL, T98, Tz} 
and some related elliptic and parabolic partial differential equations, 
\cite[\arti]{Fa99, FaHu01, HZ12, HuZ10}. 
For some recent physics applications of analysis on fractals see \cite[\arti]{ADT1,ADT2, JoeChen, FKS, KKPSS, Reuter2011, Str09}, and for analysis on fractals in general see \cite{Ba98, Ki01, Str06}. 

Once a diffusion is known to exist, we may regard its infinitesimal generator as the Laplacian $\Delta$, and employ general functional analytic tools (such as semigroup theory or variational methods, \cite{Ev98}) to solve equations of type $\Delta u=f$ and $\frac{\partial u}{\partial t}=\Delta u+F(u)$, and even such of form $\Delta \Phi(u)=f$ or $\frac{\partial u}{\partial t}=\Delta \Phi(u)$, with possibly nonlinear transformations $F$ and $\Phi$. Note that these equations do not include analogs of first order operators (gradients). However, we would like to investigate  scalar equations of type
\begin{equation}\label{E:elliptic}
\diverg(a(\nabla u))=f
\end{equation}
or 
\begin{equation}\label{E:drift}
\Delta u+b(\nabla u)=f
\end{equation}
with possibly nonlinear $a$ and $b$, or vector equations like for instance the Navier-Stokes system 
\begin{equation}\label{E:NS}
\begin{cases}
\frac{\partial u}{\partial t}+(u\cdot\nabla)u-\Delta u+\nabla p=0,\\
\diverg u=0,
\end{cases}
\end{equation}
or the magnetic Schr\"odinger equation
\begin{equation}\label{E:magnetic}
i\frac{\partial u}{\partial t}=(-i\nabla-A)^2u+Vu.
\end{equation}

Previous constructions \cite{Ki93, Ku89, St00, T00}, of first order operators related to diffusions on fractals were rather based 
on purely probabilistic and point-wise approaches, and perhaps for this reason not quite flexible enough to fit into a setup that allows to investigate partial differential equations containing first order terms. The machinery of \cite{CS03, CS09}, together with further developments in  \cite{HRT,HTa,HTb,IRT}, provides a functional analytic definition of a first order derivation (respectively gradient) and a framework suitable for a comfortable analysis of problems like (\ref{E:elliptic})-(\ref{E:magnetic}) on fractals. 

It is our aim in this paper to highlight elements of this toolkit and to announce some related results. We proceed as follows. In the next section we state our main hypotheses and collect some useful facts on Dirichlet forms and energy measures. In Section \ref{S:CS} we review the basic setup of \cite{CS03, CS09} and discuss related notions of vector analysis proposed in \cite{HRT}. First applications to scalar valued partial differential equations of types (\ref{E:elliptic}) and (\ref{E:drift}) are then presented in Section \ref{S:scalar}, and some results on analogs of (\ref{E:NS}) in Section \ref{S:vector}. In Section \ref{S:magnetic} we discuss an apporach to (\ref{E:magnetic}). We also present the definition of related Dirac operators proposed in  abstract form in \cite{CS03} and in pointwise form in\cite{HTb}.

\section{Dirichlet forms and energy measures}\label{S:Dforms}

Let $X$ be a locally compact separable metric space and $m$ a Radon measure on $X$ such that each nonempty open set is charged positively. We assume that $(\mathcal{E},\mathcal{F})$ is a symmetric local regular Dirichlet form on $L_2(X,m)$ with core $\mathcal{C}:=\mathcal{F}\cap C_0(X)$.
Endowed with the norm $\left\|f\right\|_{\mathcal{C}}:=\mathcal{E}(f)^{1/2}+\sup_X|f|$ the space $\mathcal{C}$ becomes an algebra and in particular,
\begin{equation}\label{E:boundmult}
\mathcal{E}(fg)^{1/2}\leq \left\|f\right\|_{\mathcal{C}}\left\|g\right\|_{\mathcal{C}}, \ \ f,g\in\mathcal{C},
\end{equation}
see \cite{BH91}. For any $g,h\in\mathcal{C}$ we can define a finite signed Radon measure $\Gamma(g,h)$ on $X$ such that
\[2\int_X f\:d\Gamma(g,h)=\mathcal{E}(fg,h)+\mathcal{E}(fh,g)-\mathcal{E}(gh,f)\ , \ \ f\in \mathcal{C},\]
the \emph{mutual energy measure} of $g$ and $h$. By approximation we can also define the mutual energy measure $\Gamma(g,h)$ for general $g, h\in\mathcal{F}$. Note that $\Gamma$ is symmetric and bilinear, and $\Gamma(g)\geq 0$, $g\in \mathcal{F}$. For details we refer the reader to \cite{FOT94}. We provide some examples.

\

\begin{examples}\label{Ex:Dforms}\mbox{}
\begin{enumerate}
\item[(i)] \emph{Dirichlet forms on Euclidean domains}. Let $X=\Omega$ be a bounded domain in $\mathbb{R}^n$ with smooth boundary $\partial\Omega$ and
\[\mathcal{E}(f,g)=\int_\Omega \nabla f\nabla g\:dx, \ \ f,g\in C^\infty(\Omega).\]
If $H_0^1(\Omega)$ denotes the closure of $C^\infty(\Omega)$ with respect to the scalar product $\mathcal{E}_1(f,g):=\mathcal{E}(f,g)+\left\langle f,g\right\rangle_{L_2(\Omega)}$, then $(\mathcal{E},H_0^1(\Omega))$ is a local regular Dirichlet form on $L_2(\Omega)$. The mutual energy measure of $f,g\in H_0^1(\Omega)$ is given by $\nabla f\nabla g dx$.
\item[(ii)] \emph{Dirichlet forms on Riemannian manifolds}. Let $X=M$ be a smooth compact Riemannian manifold and 
\[\mathcal{E}(f,g)=\int_M \left\langle df, dg\right\rangle_{T^\ast M}\:dvol, \ \ f,g\in C^\infty(M).\]
Here $dvol$ denotes the Riemannian volume measure. Similarly as in (i) the closure of $\mathcal{E}$ in $L_2(M,dvol)$ yields a local regular Dirichlet form. The mutual energy measure of two energy finite functions $f,g$ is given by $\left\langle df, dg\right\rangle_{T^\ast M}\:dvol$.
\item[(iii)] \emph{Dirichlet forms induced by resistance forms on fractals}. Let $X$ be a set and $(\mathcal{E},\overline{\mathcal{F}})$ a local resistance form on it such that $X$, endowed with the corresponding resistance metric $R$, is complete, separable and locally compact. For any Borel regular measure $m$ on $(X,R)$ such that $0<m(B(x,r))<\infty$, the space $(\overline{\mathcal{F}}\cap L_2(X,m), \mathcal{E}_1)$ is Hilbert, and denoting by $\mathcal{F}$ the closure of $C_0(X)\cap\overline{\mathcal{F}}$ in it, we obtain a local regular Dirichlet form $(\mathcal{E},\mathcal{F})$ on $L_2(X,m)$ (see for instance \cite[Section 9]{Ki12}). Here we have again used the standard notation $\mathcal{E}_1(f,g)=\mathcal{E}(f,g)+\left\langle f,g\right\rangle_{L_2(X,m)}$.
\end{enumerate}
\end{examples}

\begin{remark}\label{R:Kusuoka}
In Examples \ref{Ex:Dforms} (i) and (ii) the energy measures have been absolutely continuous with respect to the given reference measure. For diffusions on self-similar fractals this is typically not true if we choose the corresponding self-similar Hausdorff type measure as reference measure, see for instance \cite{BBST99} or \cite{Hino05}. We may, however, use Kusuoka type measures as reference measures to produce absolute continuity, see for instance \cite{HRT, Ka12, Ki08, Ku89, T08}.
\end{remark}

\section{$1$-forms and vector fields}\label{S:CS}

Following \cite{CS03, CS09} we consider $\mathcal{C}\otimes\mathcal{B}_b(X)$, where $\mathcal{B}_b(X)$ denotes the space of bounded Borel functions on $X$. We endow this tensor product with the symmetric bilinear form 
\begin{equation}\label{E:scalarprodH}
\left\langle a\otimes b, c\otimes d\right\rangle_{\mathcal{H}}:=\int_X bd\:d\Gamma(a,c),
\end{equation}
$a\otimes b, c\otimes d \in \mathcal{C}\otimes \mathcal{B}_b(X)$, let $\left\|\cdot\right\|_\mathcal{H}$ denote the associated seminorm on $\mathcal{C}\otimes\mathcal{B}_b(X)$ and write
\[ker\:\left\|\cdot\right\|_\mathcal{H}:=\left\lbrace \sum_i a_i\otimes b_i\in\mathcal{C}\otimes\mathcal{B}_b(X): \left\|\sum_i a_i\otimes b_i\right\|_\mathcal{H}=0\right\rbrace \]
(with finite linear combinations). To the Hilbert space $\mathcal{H}$ obtained as the completion of $\mathcal{C}\otimes\mathcal{B}_b(X)/ker\:\left\|\cdot \right\|_{\mathcal{H}}$ with respect to $\left\|\cdot\right\|_\mathcal{H}$ we refer as the \emph{space of differential $1$-forms on $X$}, cf. \cite{CS03, CS09, H11, IRT}.

The space $\mathcal{H}$ becomes a bimodule if we declare the algebras $\mathcal{C}$ and $\mathcal{B}_b(X)$ to act on it as follows: For $a\otimes b \in  \mathcal{C}\otimes \mathcal{B}_b(X)$, $c\in \mathcal{C}$ and $d\in\mathcal{B}_b(X)$ set
\begin{equation}\label{E:left}
c(a\otimes b):=(ca)\otimes b - c\otimes (ab)\ 
\end{equation}
and 
\begin{equation}\label{E:right}
(a\otimes b)d:=a\otimes (bd).
\end{equation} 
In \cite{CS03} and \cite{IRT} it has been shown that (\ref{E:left}) and (\ref{E:right}) extend to well defined left and right actions of the algebras $\mathcal{C}$ and $\mathcal{B}_b(X)$ on $\mathcal{H}$. From (\ref{E:scalarprodH}) and the Leibniz rule for energy measures, see \cite[Theorem 3.2.2]{FOT94}, it can be seen that left and right multiplication agree for any $c\in\mathcal{C}$, and as
\[\max\left\lbrace \left\|(a\otimes b)c\right\|_\mathcal{H},  \left\|c(a\otimes b)\right\|_\mathcal{H}\right\rbrace \leq \sup_X|c|\left\|a\otimes b\right\|_\mathcal{H},\] 
it follows by approximation that they agree for all $c\in\mathcal{B}_b(X)$, see \cite{IRT}. 

A \emph{derivation operator} $\partial: \mathcal{C}\to \mathcal{H}$ can be defined by setting
\[\partial f:= f\otimes \mathbf{1}.\]
It obeys the Leibniz rule,
\begin{equation}\label{E:Leibniz}
\partial(fg)=f\partial g + g\partial f, \ \ f,g \in \mathcal{C},
\end{equation}
and is a bounded linear operator satisfying
\begin{equation}\label{E:normandenergy}
\left\|\partial f\right\|_\mathcal{H}^2=\mathcal{E}(f), \ \ f\in\mathcal{C}.
\end{equation}
On Euclidean domains and on smooth manifolds the operator $\partial$ coincides with the classical exterior derivative (in the sense of $L_2$-differential forms). Details can be found in \cite{CS03, CS09, H11, HRT, IRT}.

Being Hilbert, $\mathcal{H}$ is self-dual. We therefore regard $1$-forms also as \emph{vector fields} and $\partial$ as the \emph{gradient operator}. Let $\mathcal{C}^\ast$ denote the dual space of $\mathcal{C}$, normed by
\[\left\|w\right\|_{\mathcal{C}^\ast}=\sup\left\lbrace |w(f)|: f\in\mathcal{C}, \left\|f\right\|_{\mathcal{C}}\leq 1\right\rbrace. \]
Given $f,g\in\mathcal{C}$, consider the functional
\[u\mapsto \partial^\ast(g\partial f)(u):=-\left\langle \partial u, g\partial f\right\rangle_\mathcal{H}=-\int_Xg\:d\Gamma(u,f)\]
on $\mathcal{C}$. It defines an element $\partial^\ast(g\partial f)$ of $\mathcal{C}^\ast$, to which we refer as the \emph{divergence of the vector field $g\partial f$}.  

\begin{lemma}
The divergence operator $\partial^\ast$ extends continuously to a bounded linear operator from $\mathcal{H}$ into $\mathcal{C}^\ast$
with $\left\|\partial^\ast v\right\|_{\mathcal{C}^\ast}\leq \left\|v\right\|_{\mathcal{H}}$, $v\in\mathcal{H}$. We have
\[\partial^\ast v(u)=-\left\langle \partial u, v\right\rangle_\mathcal{H}\]
for any $u\in\mathcal{C}$ and any $v\in\mathcal{H}$. 
\end{lemma}

The Euclidean identity
\[div\:(g \: grad\:f)= g\Delta f +\nabla f \nabla g\]
has a counterpart in terms of $\partial$ and $\partial^\ast$. Let $(A, dom\:A)$ denote the infinitesimal $L_2(X,\mu)$-generator of $(\mathcal{E},\mathcal{F})$. 

\begin{lemma}
We have
\[\partial^\ast (g\partial f)= g A f + \Gamma(f,g)\ ,\]
for any simple vector field $g\partial f$, $f,g\in\mathcal{C}$, and in particular, $A f=\partial^\ast \partial f$ for $f\in\mathcal{C}$.
\end{lemma}

Proofs of these results are given in \cite[Section 3]{HRT}. This \emph{distributional perspective} can be complemented by the following point of view. The operator $\partial$, equipped with the domain $\mathcal{C}$, may be seen as densely defined unbounded operator 
\[\partial: L_2(X,m)\to\mathcal{H}.\]
Since $(\mathcal{E},\mathcal{F})$ is a Dirichlet form, $\partial$ extends uniquely to a closed linear operator $\partial$ with domain $dom\:\partial=\mathcal{F}$. The divergence $\partial^\ast$, seen as an operator
\[\partial^\ast: \mathcal{H}\to L_2(X,m),\]
will be unbounded, note that in general the inclusions $\mathcal{C}\subset L_2(X,m)\subset \mathcal{C}^\ast$ are proper. As usual $v\in\mathcal{H}$ is said to be a member of $dom\:\partial^\ast$ if
there exists some $v^\ast\in L_2(X,m)$ such that $\left\langle u,  v^\ast\right\rangle_{L_2(X,m)}=-\left\langle \partial u, v\right\rangle_\mathcal{H}$ for all $u\in\mathcal{C}$. In this case $\partial^\ast v:=v^\ast$ and 
\[\left\langle u, \partial^\ast v\right\rangle_{L_2(X,m)}=-\left\langle \partial u, v\right\rangle_\mathcal{H}\ , u\in \mathcal{C},\]
i.e. $-\partial^\ast$ is the adjoint operator of $\partial$. It is immediate that $\left\lbrace \partial f: f\in dom\:A\right\rbrace\subset dom\:\partial^\ast$. As $\partial$ is densely defined and closed, the domain $dom\:\partial^\ast$ of $\partial^\ast$ is automatically dense in $\mathcal{H}$. 

We say that $(\mathcal{E},\mathcal{F})$ admits a \emph{spectral gap} if there exists some $c>0$ such that
\begin{equation}\label{E:specgap}
\int_X(f-f_X)^2dm\leq c\:\mathcal{E}(f)
\end{equation}
for any $f\in\mathcal{F}$, where $f_X=\frac{1}{m(X)}\int_Xf\:dm$. If $(\mathcal{E},\mathcal{F})$ has a spectral gap, then the image $Im\:\partial$ of $\partial$ is a closed subspace of $\mathcal{H}$. In this case the space $\mathcal{H}$ decomposes orthogonally into $Im\:\partial$ and its complement $(Im\:\partial)^\bot$, what  implies $(Im\:\partial)^\bot=ker\:\partial^\ast$, and as a consequence we observe the following explicit description of $dom\:\partial^\ast$.

\begin{corollary} Assume that $(\mathcal{E},\mathcal{F})$ admits a spectral gap, (\ref{E:specgap}). Then the domain $dom\:\partial^\ast$ agrees with 
\[\left\lbrace v\in\mathcal{H}: v=\partial f+w: f\in dom\:A\ ,\ w\in ker\:\partial^\ast\right\rbrace.\]
For any $v=\partial f+w$ with $f\in dom\:A$ and $w\in ker\:\partial^\ast$ we have $\partial^\ast v=A f$. 
\end{corollary}
The proof is short and straightforward, see \cite[Corollary 2.2]{HTb}.

\section{Scalar PDE involving first order terms}\label{S:scalar}

The results of the preceding section may be used to obtain some results on equations of type (\ref{E:elliptic}) and (\ref{E:drift}). We quote from \cite[Section 4]{HRT}. First consider the quasilinear equation
\begin{equation}\label{E:quasilinear}
\partial^\ast a(\partial u)=f.
\end{equation}
on $L_2(X,m)$. In the situation of Example \ref{Ex:Dforms} (i) it agrees with (\ref{E:elliptic}). Assume that $a:\mathcal{H}\to\mathcal{H}$ satisfies the following monotonicity, growth and coercivity conditions: 
\begin{equation}\label{E:monotone}
\left\langle a(v)-a(w),v-w\right\rangle_\mathcal{H}\geq 0\ \ \text{for all $v,w\in Im\:\partial$},
\end{equation}
\begin{equation}\label{E:growth}
\left\|a(v)\right\|_{\mathcal{H}}\leq c_0(1+\left\|v\right\|_\mathcal{H})\ \ \ \text{for all $v\in Im\:\partial$}
\end{equation}
with some constant $c_0>0$, and
\begin{equation}\label{E:coercive}
\left\langle a(v),v\right\rangle_\mathcal{H}\geq c_1\left\|v\right\|_\mathcal{H}^2-c_2 \ \ \text{for all $v\in Im\:\partial$}
\end{equation}
with constants $c_1>0$, $c_2\geq 0$. For simplicity we assume the validity of a \emph{Poincar\'e inequality},
\begin{equation}\label{E:poincare}
\left\|f\right\|_{L_2(X,m)}^2\leq c_P\:\mathcal{E}(f)
\end{equation}
with some constant $c_P>0$ for all $f\in L_2(X,m)$ with $\int_Xfdm=0$. A function $u\in \mathcal{F}$ is called a \emph{weak solution to (\ref{E:quasilinear})} if 
\[\left\langle a(\partial u),\partial v\right\rangle_\mathcal{H}=-\left\langle f,v\right\rangle_{L_2(X,m)}\ \ \text{ for all $v\in \mathcal{F}$}.\]

By classical methods, \cite[Section 9.1]{Ev98}, we obtain the following result.
\begin{theorem}
Assume $a$ satisfies (\ref{E:monotone}), (\ref{E:growth}) and (\ref{E:coercive}) and suppose (\ref{E:poincare}) holds.
Then (\ref{E:quasilinear}) has a weak solution.  
Moreover, if $a$ is strictly monotone, i.e.
\begin{equation}\label{E:strictlymon}
\left\langle a(v)-a(w),v-w\right\rangle_\mathcal{H}\geq c_3\left\|v-w\right\|_\mathcal{H}^2 \ \ \text{ for all $v,w\in Im\:\partial$}
\end{equation}
with some constant $c_3>0$, then (\ref{E:quasilinear}) has a unique weak solution.
\end{theorem}

An analog of (\ref{E:drift}) can be treated in a similar manner. Consider
\begin{equation}\label{E:nondiv}
-Au+ b(\partial u)+\varrho u=0,
\end{equation}
where $\varrho>0$ and $b$ is a generally non-linear function-valued mapping on $\mathcal{H}$. Assume that $b:\mathcal{H}\to L_2(X,m)$ is such that 
\begin{equation}\label{E:nondivgrowth}
\left\|b(v)\right\|_{L_2(X,m)}\leq c_5(1+\left\|v\right\|_\mathcal{H}),\ v\in Im\:\partial,
\end{equation}
with some $c_5>0$. A function $u\in \mathcal{F}$ is called a weak solution to (\ref{E:nondiv}) if 
\[\mathcal{E}(u,v)+\left\langle b(\partial u),\partial v\right\rangle_\mathcal{H}+\varrho\left\langle u,v\right\rangle_{L_2(X,m)}=0 \ \text{ for all $v\in \mathcal{F}$.}\] 
From \cite[Section 9.2.2, Example 2]{Ev98}, we then obtain the following.

\begin{theorem}
Assume that the embedding $\mathcal{F}\subset L_2(X,m)$ is compact and that (\ref{E:nondivgrowth}) holds. Then for any sufficiently large $\varrho>0$ there exists a weak solution to (\ref{E:nondiv}).
\end{theorem}

\section{Navier-Stokes equations}\label{S:vector}

In this section we comment on equations of type (\ref{E:NS}) which provide some more interesting applications for the notions discussed in Section \ref{S:CS}.

Together with suitable boundary conditions the Navier-Stokes system (\ref{E:NS}) describes the flow of an incompressible and homogeneous fluid in a Euclidean domain with velocity field $u$ and subject to the pressure $p$. In a one-dimensional situation it reduces to an Euler equation $\partial u/\partial t+\partial p/\partial x=0$ that has only stationary solutions. In \cite{HTa} we have proposed to investigate an analog of (\ref{E:NS}) on compact connected topologically one-dimensional fractals $X$. We collect some items necessary to formulate it.

Assume that the space $X$ is \emph{compact, connected and topologically one-dimensional} and that $(\mathcal{E},\mathcal{F})$ admits a \emph{spectral gap}, (\ref{E:specgap}). Combined with several results on Hodge decompositions and topology, cf. \cite[Sections 4,5 and 6]{HTa}, the assumption of topological one-dimensionality had motivated to define a Laplacian $\Delta_1$ on $1$-forms by 
\begin{equation}\label{E:formlaplace}
\Delta_1:=\partial\partial^\ast,
\end{equation}
seen as an unbounded operator on $\mathcal{H}$ with domain $dom\:\Delta_1=\left\lbrace \omega \in dom\:\partial^ \ast: \partial^ \ast\omega\in \mathcal{F}\right\rbrace$.

\begin{theorem}\label{T:formlaplace}
The operator $(\Delta_1, dom\:\Delta_1)$ is a self-adjoint operator on $\mathcal{H}$.
\end{theorem}

A proof can be found in \cite[Section 6]{HTb}. Theorem \ref{T:formlaplace} allows to talk about harmonic forms: A $1$-form $\omega\in\mathcal{H}$ is called \emph{harmonic} if $\omega\in dom\:\Delta_1$ and $\Delta_1\omega=0$. From compactness and topological one-dimensionality can deduce the following, cf. \cite[Theorem 6.2]{HTb}.

\begin{theorem}\label{T:harmonic}
A $1$-form $\omega\in\mathcal{H}$ is harmonic if and only if it is in $(Im\:\partial)^\bot$.
\end{theorem} 

The proof of of Theorem \ref{T:harmonic} is rather subtle, it involves a description of $(Im\:\partial)^\bot$ in terms of locally harmonic forms. We refer the reader to \cite{HTa}. Note that Theorem \ref{T:harmonic} indicates that in this situation the definition (\ref{E:formlaplace}) is appropriate.

Also for the convection term $(u\cdot\nabla)u$ in (\ref{E:NS}) we propose a substitute which by one-dimensionality seems reasonable. Our choice is motivated by the Euclidean situation: Given a vector field $u$, the quantity
\begin{equation}\label{E:weakdiv}
-\int |u|^2 \diverg v\:dx,
\end{equation}
seen as a functional on a space of test vector fields $v$, provides a formulation of $\nabla|u|^2$ in the weak sense. In our situation we set
\[dom_c\partial^\ast:=\left\lbrace v\in dom\:\partial^\ast: \partial^\ast v \in C(X)\right\rbrace\]
and given $u\in\mathcal{H}$, define
\begin{equation}\label{E:convection}
\partial\Gamma_\mathcal{H}(u)(v):=-\left\langle (\partial^\ast v) u, u\right\rangle_\mathcal{H}, \ \ v\in dom_c\partial^\ast.
\end{equation}
This seems reasonable by a \emph{fiberwise (respectively $m$-a.e. pointwise)} representation for $\mathcal{H}$ proved in \cite[Section 2]{HRT} and \cite[Theorem 2.2]{HTb}: 

\begin{theorem}\label{T:fibers} 
Let $\nu$ be a Radon measure such that all energy measures are absolutely continuous with respect to $\nu$. There are a family of Hilbert spaces $\left\lbrace \mathcal{H}_x\right\rbrace_{x\in X}$ and surjective linear maps
$\omega\mapsto \omega_x$ from $\mathcal{H}$ onto $\mathcal{H}_x$ such that the direct integral
$\int^\oplus_K \mathcal{H}_x\nu(dx)$ is isometrically isomorphic to $\mathcal{H}$ and in particular,
\[\left\|\omega\right\|_{\mathcal{H}}^2=\int_K\left\|\omega_x\right\|_{\mathcal{H},x}^2\nu(dx), \ \ \omega\in\mathcal{H}.\]
\end{theorem}

Theorem \ref{T:fibers} itself is more general, it does neither require $X$ to be compact or topologically one-dimensional nor $(\mathcal{E},\mathcal{F})$ to admit a \emph{spectral gap}.

If we replace the Euclidean norm $|\cdot|$ in (\ref{E:weakdiv}) by the norms $\left\|\cdot\right\|_{\mathcal{H},x}$ of the fibers $\mathcal{H}_x$,
we arrive at (\ref{E:convection}). On the other hand, we have the classical identity
\[\frac{1}{2}\nabla|u|^2=(u\cdot \nabla)u+u\times \curl u.\]
In a one-dimensional situation there should be no nonzero $2$-forms, hence $\curl u$ should be trivial, so that (\ref{E:convection}) is a good substitute for $(u\cdot\nabla)u$. See \cite{HTa} for a more detailed discussion. Altogether this gives a strong heuristic motivation to regard
\begin{equation}\label{E:NS2}
\begin{cases}
\frac{\partial u}{\partial t}+\frac{1}{2}\partial \Gamma_\mathcal{H}(u)-\Delta_1 u+\partial p=0\\
\partial^\ast u=0.
\end{cases}
\end{equation}
as a suitable analog of a (\ref{E:NS}) on a compact topologically one-dimensional space. Note that this is a boundary free formulation. We say that a square integrable $dom\:\partial^\ast$-valued function $u$ on $[0,\infty)$ provides a \emph{weak solution to (\ref{E:NS2})} with initial condition $u_0\in ker\:\partial^\ast$ if 
\begin{equation}\label{E:weaksol}
\begin{cases}
\left\langle u(t),v\right\rangle_\mathcal{H}-\left\langle u_0,v\right\rangle_\mathcal{H} +\int_0^t \partial\Gamma_\mathcal{H}(u(s))(v)ds+\int_0^t\left\langle\partial^\ast u(s),\partial^\ast v\right\rangle_{L_2(X,m)}ds=0\\
\partial^\ast u(t)=0
\end{cases}
\end{equation}
for a.e. $t\in [0,\infty)$ and all $v\in ker\:\partial^\ast$. By some immediate simplifications we then observe stationarity and uniqueness of solutions. In other words, the behavior of the system on a compact topologically one-dimensional space resembles its behavior on a one-dimensional Euclidean space. 

\begin{theorem}\mbox{}
Any weak solution $u$ of (\ref{E:NS2}) is harmonic and stationary, i.e. $u$ is independent of $t\in [0,\infty)$. Given an initial condition $u_0$  the corresponding weak solution is uniquely determined.
\end{theorem}

Note that we have not made any restriction on the Hausdorff dimension $d_H$ of $X$. Indeed there are examples of spaces of \emph{any Hausdorff dimension $1\leq d_H<\infty$} such that the previous Theorem holds. It is the topological dimension that governs the behavior of (\ref{E:NS2}).

\begin{remark}\label{R:fibers}
Logically Theorem \ref{T:fibers} is not needed to set up the model (\ref{E:NS2}), we have just included it here to support the intuition behind our choice of substitute terms. We would also like to remark that even though the energy measures might not be absolutely continuous with respect to the initial reference measure $m$, one can always construct a finite Radon measure $\nu$ that has this property.
\end{remark}

For the rest of this section we specialize further to the situation of Examples \ref{Ex:Dforms} (iii), that is, we assume $X$ to be a set and $(\mathcal{E},\overline{\mathcal{F}})$ to be a resistance form on it such that $X$, together with the resistance metric $R$, is a compact and connected metric space. We further assume that $m$ is a Borel regular measure on $X$ as in Examples \ref{Ex:Dforms} (iii) so that consequently a local regular Dirichlet form $(\mathcal{E},\mathcal{F})$ is obtained by taking the closure. Then all our previous results may be applied for $(\mathcal{E},\mathcal{F})$. We finally assume that $(X,R)$ is topologically one-dimensional.

\begin{remark}
We conjecture that any set that carries a regular resistance form is a topologically one-dimensional space when equipped with the associated resistance metric.
\end{remark}

Note that in the resistance form case points have positive capacity. This property allows to prove the following equivalence, \cite[Section 5]{HTa}.

\begin{theorem}
Assume that points have positive capacity and topological dimension is one. Then a nontrivial solution to (\ref{E:NS2}) exists if and only if the first \v{C}ech cohomology $\check{H}^1(X)$ of $X$  is nontrivial.
\end{theorem}

In the resistance form context Neumann derivatives are well-defined, and it is not difficult to see that if the Navier-Stokes system (\ref{E:NS2}) is considered with a nonempty boundary, it may have additional nontrivial solutions arising from solutions of a related Neumann problem.   

Let $B\subset X$ be a finite set, which is interpreted as the boundary of $X$. By $G_B$ we denote the Green operator associated with the boundary $B$ with respect to $(\mathcal{E},\overline{\mathcal{F}})$, \cite[Definition 5.6]{Ki03}, and $\mathcal{D}_{B,0}^L$ its image in $\overline{\mathcal{F}}$. Let $\mathcal{H}_B$ denote the $B$-harmonic functions with respect to $(\mathcal{E},\overline{\mathcal{F}})$, \cite[Definition 2.16]{Ki03}, and note that $\overline{\mathcal{F}}=\overline{\mathcal{F}}_B\oplus\mathcal{H}_B$, where
\[\overline{\mathcal{F}}_B:=\left\lbrace u\in\overline{\mathcal{F}}: u|_B=0\right\rbrace.\]
A $B$-harmonic function $h$ is harmonic on $B^c$ in the Dirichlet form sense, more precisely, it satisfies
$\mathcal{E}(h,\psi)=0$
for all $\psi\in\overline{\mathcal{F}}_B$. The space $\mathcal{D}^L:=\mathcal{D}^L_{B,0}+\mathcal{H}_B$ is seen to be independent of the choice of $B$, \cite[Theorem 5.10]{Ki03}. For any $u\in \mathcal{D}^L$ and any $p\in X$ the Neumann derivative $(du)_p$ of $u$ at $p$ can be defined, \cite[Theorems 6.6 and 6.8]{Ki03}. If $\varphi$ is a function on $B$, then a function $h_\varphi\in\overline{\mathcal{F}}$ is called a {solution to the Neumann problem on $B^c$ with boundary values $\varphi$} if it 
is harmonic on $B^c$ and satisfies 
$(dh)_p=\varphi(p)$
for all $p\in B$. 
Such a Neumann solution $h_\varphi$ exists and is unique if and only if $\varphi$ is such that $$\sum_{p\in B}\varphi(p)=0.$$ 

We use the notation $\mathcal{H}(B^c)=\clos\lin\left\lbrace v\mathbf{1}_{B^c}:v\in\mathcal{H}\right\rbrace$. A square integrable $dom\:\partial^\ast$-valued function $u$ on $[0,\infty)$ provides a \emph{weak solution to (\ref{E:NS2}) on $B^c$} if 
\begin{equation}\label{E:NSboundary}
\begin{cases}
\left\langle u(t),v\right\rangle_\mathcal{H}-\left\langle u(0),v\right\rangle_\mathcal{H} +\int_0^t \partial\Gamma_\mathcal{H}(u(s))(v)ds+\int_0^t\left\langle \partial^\ast u(s),\partial^\ast v\right\rangle_{L_2(X,m)}ds=0\\
\left\langle u(t), \partial\psi\right\rangle_\mathcal{H}=0
\end{cases}
\end{equation}
for a.e. $t\in [0,\infty)$, all $v\in \text{dom}\:\partial^\ast \cap \mathcal{H}(B^c)$ and all $\psi\in\overline{\mathcal{F}}_B$.

\begin{theorem}
Assume that points have positive capacity and topological dimension is one.  If $h$ is the unique, up to an additive constant, harmonic function on $B^c$ with normal derivatives $\varphi$ on $B$, then 
\[u(t)=\partial h,\ \  t\in [0,\infty).\]
is the unique weak solution to (\ref{E:NS2}) on $B^c$ with the Neumann boundary values $\varphi$ on~$B$.
\end{theorem}

\begin{remark}
In (\ref{E:NSboundary}) we have considered weak solutions to (\ref{E:NS2}). For weak solutions the pressure $p$ does not occur explicitely. However, any definition of strong solution to \ref{E:NS} should lead to the relation
\[p(t)=-\frac12\Gamma(h),\ \  t\in[0,\infty),\]
seen as an equality of measures. A similar statement could be written for the boundary free case (\ref{E:weaksol}).
\end{remark}

For more details and the Hodge theory leading to the statements of this section we refer the reader to \cite{HTb}.

\section{Magnetic Schr\"odinger equations}\label{S:magnetic}

We turn to results concerning the magnetic Schr\"odinger equation (\ref{E:magnetic}). To discuss this equation we do not need to assume that $X$ is compact or topologically one-dimensional. As in Section \ref{S:Dforms} it may just be an arbitrary locally compact separable metric space equipped with a Radon measure $m$ that charges any nonempty open set positively and carrying a symmetric local regular Dirichlet form $(\mathcal{E},\mathcal{F})$ on $L_2(X,m)$ with core $\mathcal{C}:=\mathcal{F}\cap C_0(X)$. However, to investigate (\ref{E:magnetic}) we will now assume that $(\mathcal{E},\mathcal{F})$ \emph{possesses energy densities with respect to the reference measure $m$}, i.e. for any $g,h\in \mathcal{F}$ the measure $\Gamma(g,h)$ is absolutely continuous with respect to $m$. 

\begin{remark}
As previously mentioned in Remark \ref{R:fibers} we can always construct a measure $\nu$ with respect to which all energy measures are absolutely continuous. For the cases that the given Dirichlet form $(\mathcal{E},\mathcal{F})$ on $L_2(X,m)$ is transient or induced by a resistance form we have shown in \cite{HRT} that $(\mathcal{E},\mathcal{C})$ is a closable form on the space $L_2(X,\nu)$ of functions that are square integrable with respect to this new measure $\nu$. Then this change of measure merely amounts to a change of domains. It is not difficult to show that under mild assumptions $(\mathcal{E},\mathcal{C})$ is always closable with respect to this measure $\nu$. We will discuss this matter in a later paper.
\end{remark}

In \cite{HTb} we have studied analogs of the magnetic Hamiltonian $(-i\nabla-A)^2+V$ and in particular, have verified their essential self-adjointness. To sketch this result, let $L_{2,\mathbb{C}}(X,m)$, $\mathcal{F}_\mathbb{C}$, $\mathcal{C}_\mathbb{C}$ and $\mathcal{H}_\mathbb{C}$ denote the natural complexifications of $L_2(X,m)$, $\mathcal{F}$, $\mathcal{C}$ and $\mathcal{H}$, respectively. The natural extensions to the complex case of $\mathcal{E}$ and the corresponding energy measures $\Gamma(f,g)$ are again denoted by the same symbol. Note that they are conjugate symmetric and linear in the first argument. If both arguments agree, they yield a real nonnegative number and a real nonnegative measure, respectively. The first result concerns related quadratic forms. Here we use the notion of quadratic form in the sense of \cite[Section VIII.6]{RS1}.

\begin{proposition}\label{P:quadratic} 
Let $a\in \mathcal{H}$ and $V\in L_\infty(X,m)$. The form $\mathcal{E}^{a,V}$, given by 
\[\mathcal{E}^{a,V}(f,g)=\left\langle (-i\partial-a)f,(-i\partial-a)g\right\rangle_{\mathcal{H}}+\left\langle fV,g\right\rangle_{L_2(X,m)}, \ \ f,g\in\mathcal{C}_\mathbb{C},\]
defines a quadratic form on $L_{2,\mathbb{C}}(X,m)$.
\end{proposition}

Proposition \ref{P:quadratic} is a slight variation of Proposition 4.1 in \cite{HTb} and up to inessential details it has the same proof. Here $a$ is seen as the \emph{magnetic vector potential} replacing $A$ in (\ref{E:magnetic}) and $V$ is the \emph{electric scalar potential}.

Now recall the fiberwise representation of $\mathcal{H}$ from Theorem \ref{T:fibers}. We define the \emph{space of real vector fields of bounded length} by
\[\mathcal{H}_\infty:=\left\lbrace v=(v_x)_{x\in X}\in\mathcal{H}: \left\|v_\cdot\right\|_{\mathcal{H},\cdot}\in L_\infty(X,m)\right\rbrace.\]
If the potential $a$ is recruited from $\mathcal{H}_\infty$ then we can obtain the closedness of $\mathcal{E}^{a,V}$ and the essential self-adjointness of the associated operator from straightforward perturbation arguments, cf. \cite[Theorem 4.1]{HTb}. Recall that $A$ denotes the generator of $(\mathcal{E},\mathcal{F})$. We denote its complexification by the same symbol.

\begin{theorem}\label{T:magnetic}
Let $a\in \mathcal{H}_\infty$ and $V\in L_\infty(X,m)$. 
\begin{enumerate}
\item[(i)] The quadratic form $(\mathcal{E}^{a,V},\mathcal{F}_\mathbb{C})$ is closed. 
\item[(ii)] The self-adjoint non-negative definite operator on $L_{2,\mathbb{C}}(X,m)$ uniquely associated with $(\mathcal{E}^{a,V},\mathcal{F}_\mathbb{C})$ 
is given by 
\[H^{a,V}=(-i\partial-a)^\ast(-i\partial-a)+V,\]
and the domain of the operator $A$ is a domain of essential self-adjointness for $H^{a,V}$. 
\end{enumerate}
\end{theorem}

The operator $H^{a,V}$ is a natural generalization of the quantum mechanical Schr\"odinger Hamiltonian $(-i\nabla-A)+V$ from (\ref{E:magnetic}). By Theorem \ref{T:magnetic} we have established a suitable framework to study a fractal counterpart 
\[i\frac{\partial u}{\partial t}=H^{a,V}u\]
of the evolution equation (\ref{E:magnetic}).

As they are closely related to magnetic Hamiltonians, we conclude this section by a brief look at Dirac operators. We have introduced a local Dirac operator in \cite{HTb}. Up to sign and complexity conventions it is defined as a matrix operator 
\begin{equation}\label{E:DiracIntro}
D=\left(
\begin{array}{ll}
0 &\partial^\ast\\
\partial &0
\end{array}\right),
\end{equation}
acting on $H_0\oplus H_1=L_2(X,m)\oplus \mathcal H$. We consider $D$ as an unbounded linear operator with domain $dom\:D:=\mathcal{F}\oplus dom\:\partial^\ast$ and have the following result, obtained in  abstract form in \cite{CS03}, and in pointwise form in\cite{HTb}. 

\begin{theorem}\label{T:Dirac}
The operator $(D, dom\:D)$ is self-adjoint operator on $H_0\oplus H_1$.
\end{theorem}

Note that as a consequence we also obtain a local matrix Laplacian $D^2$ acting on $H_0\oplus H_1$. 

According to \cite{HTa,HTb}, this Dirac operator is naturally related to the topological structure of the fractals space and, in a certain natural  a sense, to the differential geometry of the fractal (see \cite[\arti]{CI,CIL} for a discussion of the notion of a Dirac operator in the context of non-commutative analysis). In particular, for {\itshape topologically one dimensional fractals (of arbitrary Hausdorff and spectral dimensions)} our  Dirac operator gives rise  to a natural Hodge Laplacian \ $\partial \partial^*+\partial^*\partial$ \ on the appropriate differential complex.  
It will be the subject of future work to study the Hodge Laplacian for higher order differential forms defined in the probabilistic or Dirichlet form sense.  

As a side remark we note that sometimes there may be a different convention for the Dirac operator in a complex setup. For instance
\[D=\left(\begin{array}{rr}
0 & -i\partial^\ast\\
-i\partial & 0\end{array}\right)\]
has signs and imaginary factors which are somewhat more suitable in relation to the magnetic magnetic Schr\"odinger operator $H^{a,V}$. 

Finally, we would like to point out related perturbation results. Assume that $b\in\mathcal{H}_\infty$ and set
\[\mathcal{Q}(f,g):=\mathcal{E}(f,g)-\int_X g(x)\left\langle b_x,\partial_xf\right\rangle_{\mathcal{H}_x}m(dx),\]
$f,g\in\mathcal{F}$. Here $b_x$ and $\partial_xf$ denote the images of $b$ and $\partial f$ under the projection from $\mathcal{H}$ onto $\mathcal{H}_x$ as in Theorem \ref{T:fibers}. For $\alpha\geq 0$ write
\[\mathcal{Q}(f,g):=\mathcal{Q}(f,g)+\alpha\left\langle f,g\right\rangle_{L_2(X,m)}.\]
We may then conclude the following.

\begin{theorem}\mbox{}
\begin{enumerate}
\item[(i)] For any $\alpha\geq 0$ the form $(\mathcal{Q}_\alpha,\mathcal{F})$ is closed on $L_2(X,m)$. It generates a strongly continuous semigroup of bounded operators on $L_2(X,m)$.
\item[(ii)] If $\alpha>0$ is sufficiently large then the associated semigroup is positivity preserving.
\item[(iii)] The generator $L^\mathcal{Q}$ of $\mathcal{Q}$ is given by
\[L^\mathcal{Q} u(x)=Au(x)+\left\langle b_x,\partial_xf\right\rangle_{\mathcal{H}_x}, \ \ u\in dom\:A.\]
\end{enumerate}
\end{theorem}
See \cite[Section 10]{HRT} and the references therein, in particular \cite{FK04}.

\end{document}